\documentclass[legalpaper]{article}
\usepackage{amsmath}
\usepackage{amssymb}
\usepackage{amsthm}
\usepackage{graphicx} 
\usepackage{tikz}
\usepackage{ifmtarg}
\usepackage[mathscr]{eucal}

\newtheorem{Theorem}[equation]{Theorem}

\newtheorem{Lemma}[equation]{Lemma}

\theoremstyle{definition}

\newtheorem{Example}[equation]{Example}

\newtheorem{Conjecture}[equation]{Conjecture}

\theoremstyle{remark}
\newtheorem{Remark}[equation]{Remark}
\newtheorem{Remarks}[equation]{Remarks}

\newcommand{\defeq}{\overset{\operatorname{\scriptstyle def.}}{=}}
\newcommand{\CC}{{\mathbb C}}

\newcommand{\ZZ}{{\mathbb Z}}

\newcommand{\SL}{\operatorname{\rm SL}}

\newcommand{\GL}{\operatorname{GL}}
\newcommand{\PGL}{\operatorname{PGL}}
\newcommand{\SO}{\operatorname{\rm SO}}

\newcommand{\grpSp}{\operatorname{\rm Sp}}

\newcommand{\algsl}{\operatorname{\mathfrak{sl}}} 

\newcommand{\Spec}{\operatorname{Spec}\nolimits}

\newcommand{\Hom}{\operatorname{Hom}}

\newcommand{\Ker}{\operatorname{Ker}}

\newcommand{\MR}[1]{}
\newcommand{\Wedge}{{\textstyle \bigwedge}}

\newcommand{\dslash}{/\!\!/}
\newcommand{\bM}{\mathbf M}
\newcommand{\bN}{\mathbf N}

\newcommand{\tslash}{/\!\!/\!\!/}
\newcommand{\tslabar}{\mathbin{
\setbox0=\hbox{/\!\!/\!\!/}\rule[0.4\ht0]{\wd0}{.3\dp0}\kern-\wd0\box0}}

\newcommand{\Gr}{\mathrm{Gr}}

\newcommand{\cR}{\mathcal R}

\newcommand{\cT}{\mathcal T}

\newcommand{\cO}{\mathcal O}

\makeatletter
\newcommand{\cA}[1][{}]{%
  \@ifmtarg{#1}%
  {\mathcal A}
  {\mathcal A(#1)}
}
\newcommand{\cAh}[1][{}]{%
  \@ifmtarg{#1}%
  {\mathcal A_\hbar}
  {\mathcal A_\hbar(#1)}
}
\makeatother

\DeclareSymbolFont{symbolsC}{U}{pxsyc}{m}{n}
\SetSymbolFont{symbolsC}{bold}{U}{pxsyc}{bx}{n}
\DeclareFontSubstitution{U}{pxsyc}{m}{n}
\DeclareMathSymbol{\medcirc}{\mathbin}{symbolsC}{7}

\newcommand{\bG}{\mathbf G}

\tikzset{SO/.style={draw,
    minimum size=12pt,inner sep=0pt, outer sep=0pt,font={$+$}}}
\tikzset{Sp/.style={draw,
    minimum size=12pt,inner sep=0pt, outer sep=0pt,font={$-$}}}
\tikzset{GL/.style={draw,
    minimum size=12pt,inner sep=0pt, outer sep=0pt}}
\tikzset{pm/.style={draw,
    minimum size=12pt,inner sep=0pt, outer sep=0pt,font={$\pm$}}}
\tikzset{mp/.style={draw,
    minimum size=12pt,inner sep=0pt, outer sep=0pt,font={$\mp$}}}

\newcommand{\Lie}{\operatorname{Lie}}

\newcommand{\po}{\ar@{}[dr]|{\text{\pigpenfont R}}}
\newcommand{\pb}{\ar@{}[dr]|{\text{\pigpenfont J}}}
\newcommand{\pp}{\ar@{}[dr]|{\text{\pigpenfont P}}}

\newcommand{\cM}{\mathcal M}

\newcommand{\CN}{{\mathcal{N}}}

\newcommand{\cS}{{\mathcal{S}}}

\newcommand{\bv}{{\mathbf v}} 
\newcommand{\bw}{{\mathbf w}} 

\numberwithin{equation}{section}

\title{S-dual of Hamiltonian $\bG$ spaces and relative Langlands duality}
\author{Hiraku Nakajima}
\date{September 2024}

\begin{document}

\maketitle
\thispagestyle{empty}

\section{Introduction}\label{sec:intro}

Let $\bG$ be a complex reductive group, and $\bM$ be
a smooth affine algebraic symplectic manifold with a 
Hamiltonian $\bG$ action. 
Physicists associate 3-dimensional $\mathcal N=4$ supersymmetric quantum field
theory\footnote{This theory
can be coupled with $\bG$-connections, the resulted theory (gauged $\sigma$-model)
is denoted by $\mathscr T_{\bG,\bM}\tslabar\bG$, and should not be confused with
$\mathscr T_{\bG,\bM}$.
In quiver gauge theories, $\bG$ is the group for a framing
vector space, while $\mathscr T_{\bG,\bM}\tslabar\bG$ is the
gauge theory obtained by changing the framing vector space
to the usual one.} $\mathscr T_{\bG,\bM}$,
the nonlinear $\sigma$-model,
to the pair $(\bG\curvearrowright\bM)$. 
This $\mathscr T_{\bG,\bM}$ appears as a boundary condition of 4-dimensional $\mathcal N=4$
super Yang-Mills theory with gauge group $\bG$.
The 4-dimensional $\mathcal N=4$ super Yang-Mills theory has S-duality (electric-magnetic duality) with $\bG$ replaced by its Langlands dual group $\bG^\vee$.
Under the S-duality, the boundary condition $\mathscr T_{\bG,\bM}$ is mapped to a
boundary condition for the 4-dimensional super Yang-Mills theory with gauge group $\bG^\vee$.
In \cite[(3.13)]{MR2610576}, Gaiotto and Witten 
observed that the S-dual boundary condition is closely
related to 3d-mirror symmetry. In many examples, the S-dual theory is
associated with another pair $(\bG^\vee\curvearrowright\bM^\vee)$.

Recall that Kapustin-Witten \cite{MR2306566} explained that the geometric Langlands conjecture is a consequence,
in a physics level of rigor, 
of the S-duality of 4-dimensional super Yang-Mills theory.
In \cite{2016arXiv160909030G}, Gaiotto explained consequences of this S-duality of boundary conditions to
the geometric Langlands conjecture in the framework
of \cite{MR2306566}. Let us summarize the main statement. 
First claim says that $(\bG\curvearrowright\bM)$
gives rise to a pair 
$(\mathcal L(C,\bG\curvearrowright\bM),
\mathcal V(C,\bG\curvearrowright\bM))$
of objects in pairs of the categories for the geometric
Langlands conjecture for a Riemann surface $C$. Roughly,
$\mathcal L(C,\bG\curvearrowright\bM)$ is a holomorphic 
lagrangian subvariety in the moduli space of $\bG$-Higgs bundles
over $C$, and $\mathcal V(C,\bG\curvearrowright\bM)$ is
a holomorphic vector bundle over the moduli space of
flat $\bG$-bundles. Geometric Langlands conjecture predicts
an equivalence of two categories with $\bG$ replaced by $\bG^\vee$.
The main statement is the following relations among objects
under the geometric Langlands equivalence:
\begin{equation}\label{eq:7}
\begin{gathered}
    \mathcal L(C,\bG\curvearrowright\bM)\longleftrightarrow
    \mathcal V(C,\bG^\vee\curvearrowright\bM^\vee),\\
    \mathcal V(C,\bG\curvearrowright\bM)\longleftrightarrow
    \mathcal L(C,\bG^\vee\curvearrowright\bM^\vee).
\end{gathered}
\end{equation}

On the other hand, the author, together with Braverman and Finkelberg, studied
\emph{ring objects} in the equivariant derived Satake category in \cite{2017arXiv170602112B}
to gain flexibility for the construction of Coulomb branches of 3d $\mathcal N=4$ gauge
theories.
This study turned out to be an appropriate framework to approach \cite{MR2610576} in
a mathematically rigorous way for two reasons: (a) it discusses Coulomb branches, which 
are exchanged with Higgs branches of paired theories in 3d mirror symmetry. (b) it
discusses the \emph{regular sheaf}, which is supposed to correspond to 
the `kernel' theory $\cT[\bG]$ in \cite[(3.13)]{MR2610576}.
In particular, the S-dual $(\bG^\vee\curvearrowright \bM^\vee)$ of the pair
$(\bG\curvearrowright \bM)$
was implicitly introduced in \cite[\S5]{2017arXiv170602112B}
and explicitly in \cite[(1.3)]{MR4286060}
(both under the assumption $\bM=T^*\bN$).
See \eqref{eq:1}.
This definition is mathematically rigorous, but 
it should be regarded as a temporary one: 
the assumption $\bM=T^*\bN$ should be weakened.
There are results \cite{bdfrt,2022arXiv220901088T} in this direction, but under the assumption that $\bM$
is a symplectic \emph{representation} of $\bG$.
Our $\bM^\vee$ is \emph{singular} in general, and
it is not clear how to define the S-dual of $(\bG^\vee\curvearrowright \bM^\vee)$ in such a case. Even
if it is nonsingular, it is unclear why we expect this operation is a \emph{duality}, i.e., 
the S-dual of $(\bG^\vee\curvearrowright \bM^\vee)$ is the original $(\bG\curvearrowright \bM)$.

In \cite{BZSV}, Ben-Zvi, Sakellaridis and Venkatesh proposed
the relative Langlands duality conjecture.
Their proposal contains many new and deep ideas, but from the author's \emph{naive}\footnote{The author does \emph{not} work on Langlands duality.} point of view, 
the following two items are important:
\begin{itemize}
    \item There is a class of $\bG\curvearrowright\bM$ (\emph{hyperspherical varieties}), which is preserved under the S-duality. (See \cite[\S3]{BZSV}.)
    \item Relations \eqref{eq:7} should hold at all ``tiers'' of 
    the Langlands program (global, local, geometric, arithmetic,
    etc.). (See \cite[Introduction]{BZSV}.) 
\end{itemize}
The importance of the first item is clear: our
definition $\bM^\vee$ yields a singular space 
in general. Regarding the second item, note that
an equivalence of categories for the local case
\emph{recovers} the definition \eqref{eq:1}.
See \cite[\S8]{BZSV}.

We discuss examples which 
are \emph{not} hyperspherical in general in this article. 
We hope that our view point gives a step towards a generalization
of relative Langlands duality to more general situations.

\subsection*{Acknowledgement}
I thank Alexander Braverman and Michael Finkelberg for
discussion on the S-duality when we wrote the
paper \cite{2017arXiv170602112B} and also afterwards.
I also thank Eric Chen and Ivan Losev for discussion on relative Langlands
duality and hyperspherical conditions respectively.
I also thank the organizers and speakers of
Summer School on Relative Langlands Duality in Minnesota in June 2024. I explained some parts of this article there, and
learned remaining parts from talks by other speakers.

\section{Lagrangians and vector bundles}

The main focus of this article is the definition and examples
of S-duality $(\bG\curvearrowright\bM) \leftrightarrow
(\bG^\vee\curvearrowright\bM^\vee)$. Therefore we describe
$\mathcal L(C,\bG\curvearrowright\bM)$ and $\mathcal V(C,\bG\curvearrowright\bM)$ only briefly. See \cite[\S2.2]{2016arXiv160909030G} for more detail. We also omit
mathematically rigorous treatments.

Regarding $\mathcal L(C,\bG\curvearrowright\bM)$, we assume that
$\bM$ has a $\CC^\times$ action scaling the symplectic form
with weight $2$. In particular, the moment map
$\mu\colon\bM\to\operatorname{Lie}\bG^*$ is equivariant if
$\CC^\times$ acts on $\operatorname{Lie}\bG^*$ with weight $2$. We consider the moduli space of triples
\begin{itemize}
    \item $\mathcal P_{K_C^{1/2}}$ : a $\CC^\times$ principal bundle, whose associated line bundle is a square root $K_C^{1/2}$ of the canonical bundle of $C$. (We consider this as a part of data, hidden in the notation $C$.)
    
    \item $\mathcal P$ : a $\bG$ principal bundle on $C$.

    \item $Z$ : a section of the associated bundle 
    $(\mathcal P\times_C\mathcal P_{K_C^{1/2}})\times^{\bG\times\CC^\times}\bM$.
\end{itemize}
We apply the moment map $\mu$ fiberwise to $Z$ to get
a Higgs bundle $(\mathcal P,\mu(Z))$. Then
$\{ (\mathcal P,Z)\}/\text{isom}\to\{(\mathcal P,\mu(Z))\}/\text{isom}$ gives a lagrangian subvariety
$\mathcal L(C,\bG\curvearrowright\bM)$. 
See \cite[1(ii)]{2015arXiv151003908N} for
a related construction. See also \cite{MR3855670}
for a mathematically rigorous treatment. (In particular,
$\{ (\mathcal P,Z)\}/\text{isom}\to\{(\mathcal P,\mu(Z))\}/\text{isom}$ should be understood in the
framework of the derived symplectic geometry.)

\begin{Example}\label{ex:lag}
    \textup{(1)} Suppose $\bM = 0$. Then 
    $\{ (\mathcal P,Z)\}/\text{isom}$ is the moduli space of
    $\bG$-bundles, embedded into the moduli space of $\bG$-Higgs
    bundles by setting Higgs fields $=0$.

    \textup{(2)} Let $\bM$ be $\bG\times\mathcal S$, where
    $\mathcal S$ is the Kostant-Slodowy slice to a principal nilpotent orbit of $\bG$. (See the paragraph after Theorem~\ref{thm:3} for the definition.) Then
    $\{ (\mathcal P,Z)\}/\text{isom}$ is nothing but
    Hitchin section introduced in \cite[\S5]{MR1174252}.
    (The point $e$ of $\mathcal S$ is identified with
    the $\SL_2$-Higgs field $(K_C^{1/2}\oplus K_C^{-1/2}, \begin{pmatrix}
        0 & 0 \\ 1 & 0
    \end{pmatrix})$ via $\SL_2\to\bG$.)
\end{Example}

Let us turn to $\mathcal V(C,\bG\curvearrowright\bM)$. We assume
$\bM = \bN\oplus\bN^*$ with a representation $\bN$ of $\bG$.
We consider a family of flat connections of the associated
vector bundle over $C$ for the representation $\bN$ of $\bG$,
parametrized by the moduli space of flat $\bG$-bundles.
We assume that the twisted deRham cohomology vanishes
except degree $1$. Then the degree $1$ cohomology groups
form a holomorphic vector bundle $\mathcal D_\bN$ over the
moduli space, and define $\mathcal V(C,\bG\curvearrowright\bM)$
as $\det\mathcal D_\bN^{-1/2}\Wedge^*\mathcal D_\bN$.
(See \cite[(2.6),(10.1)]{2016arXiv160909030G}.)

\begin{Remark}
    In \cite[\S2.2]{2016arXiv160909030G}, it was emphasized that $\mathcal V(C,\bG\curvearrowright\bM)$ should be defined as a
    BBB brane, whose mathematically rigorous treatment seems to
    be missing. It 
    was proposed roughly as a generalization of a $C^\infty$ 
    vector bundle with a connection whose curvature
    is of type (1,1) for any complex structure in the moduli
    space of flat $\bG$-bundles, viewed as a hyper-K\"ahler manifold.
\end{Remark}

\begin{Example}\label{ex:Hitchin}
    Suppose $\bM=0$, and hence $\bN=0$. Then $\mathcal D_\bN = 0$,
    and $\mathcal V(C,\bG\curvearrowright\bM)$ is the structure
    sheaf over the moduli space of flat $\bG$-bundles. The
    S-dual of this example is supposed to be Example~\ref{ex:lag}(2)
    with $\bG$ replaced by $\bG^\vee$. 
    See Example~\ref{ex:2}(2).
    Therefore the Hitchin section is corresponds to the structure sheaf under
    \eqref{eq:7}. This is a well-known expectation in geometric
    Langlands correspondence.
\end{Example}

\section{Definition of S-dual}

\subsection{Coulomb branches}\label{subsec:Coulomb}

Let us quickly recall the definition
of Coulomb branches of 3d $\mathcal N=4$
supersymmetric gauge theories associated
with a complex reductive group $\bG$
and a $\bG$-space $\bM$. 
See an exposition \cite{2022arXiv220108386N} for further information.
In \cite{2016arXiv160103586B} we define the Coulomb branch
for $(\bG,\bM)$ under the assumption $\bM = \bN\oplus\bN^*$ with a representation
$\bN$ of $\bG$. However, the actual assumption we need is $\bM = T^*\bN$ for
a smooth affine algebraic variety $\bN$.

Let $F = \CC((z))$ and $\cO = \CC[[z]]$ be 
the field of formal Laurent power series and 
the formal power series ring respectively.

Let $\Gr_\bG$ be the affine Grassmannian of $\bG$, i.e.\ 
$\bG(F)/\bG(\cO)$. We define the \emph{variety of triples} $\cR_{\bG,\bN}$ by
\begin{equation*}
    \cR_{\bG,\bN} = \left\{ [g(z),s(z)]\in \bG(F)\times^{\bG(\cO)}\bN(\cO)\,\middle|\,
    g(z)s(z)\in\bN(\cO)\right\}.
\end{equation*}
It is equipped with a projection $\pi\colon\cR_{\bG,\bN}\to\Gr_\bG$. We have an action
of $\bG(\cO)$ by the multiplication to $g(z)$ from the left. We consider the equivariant
Borel-Moore homology group $H_*^{\bG(\cO)}(\cR_{\bG,\bN})$.
It is equipped with the convolution product, which is commutative.
(The definition of the convolution product uses the smoothness of
$\bN$, hence it is not \emph{clear} how to generalize the definition
to singular $\bN$.)
We define the Coulomb branch by
\begin{equation*}
    \cM_C(\bG,\bM) \defeq \Spec H_*^{\bG(\cO)}(\cR_{\bG,\bN}).
\end{equation*}

\begin{Remark}\label{rem:just}
    We need a justification for the notation $\cM_C(\bG,\bM)$,
    as the above definition uses $\bN$, rather than $M = T^*\bN$. At least when $\bN$ is a representation of $\bG$, it was proved that the result depends on $\bM = \bN\oplus\bN^*$ up to
    an isomorphism in \cite[\S6(viii)]{2016arXiv160103586B}.
\end{Remark}

\begin{Example}\label{ex:1}
(1) Suppose $\bG = \CC^\times$, $\bM = 0$ (hence $\bN=0$). Then $\cR_{\bG,\bN} = \Gr_{\bG}$.
It is well-known that $\Gr_{\bG}\cong \{ [z^n] \mid n\in\ZZ\}$. Thus 
\begin{equation*}
    H_*^{\bG(\cO)}(\cR_{\bG,\bN}) = \bigoplus_{n\in\ZZ} H_*^{\bG}(\mathrm{pt}) = 
    \bigoplus_{n\in\ZZ} \CC[w]r_n,
\end{equation*}
where $r_n$ is the fundamental class of the point $[z^n]$, and wee have used
$H_*^{\bG}(\mathrm{pt}) \cong\CC[w]$. One can compute $r_m r_n = r_{m+n}$, hence
$H_*^{\bG(\cO)}(\cR_{\bG,\bN}) = \CC[w,r_1,r_{-1}]/(r_1 r_{-1} = 1)$. Thus
$\cM_C$ is $\CC\times\CC^* = T^*(\CC^\times)$.

(2) Suppose $\bG = \CC^\times$, $\bM = \CC^2 = T^*\CC$ (hence $\bN=\CC$ and $\bN(\cO) = \CC[[z]]$). The $\bG$-action on $\bN$ is the standard weight $1$ representation. Then
\begin{equation*}
    \cR_{\bG,\bN} = \bigsqcup_{n\in\ZZ} \{ \left([z^n],s'(z)\right)
    \in\Gr_\bG\times \CC[[z]]\mid z^{-n}s'(z)\in\CC[[z]] \},
\end{equation*}
where $s'(z) = g(z)s(z)$ in the definition above. If $n \le 0$, the condition
$z^{-n}s'(z)\in\CC[[z]]$ is automatically satisfied. Therefore the fiber is $\CC[[z]]$.
On the other hand, the fiber is $z^n\CC[[z]]$ if $n > 0$. Since $z^n\CC[[z]]$ is a codimension
$n$ subspace in $\CC[[z]]$, the fundamental class of the fiber is $r'_n = w^n r_n$ compared with
the example above. (We interpret $r_n$ as the fundamental class of $\CC[[z]]$ by Thom isomorphism.) On the other hand, $r'_n = r_n$ for $n\le 0$. Therefore we have $w^1 r_1 \cdot r_{-1} = w$ in this case, hence
$H_*^{\bG(\cO)}(\cR_{\bG,\bN}) = \CC[w,r'_1,r'_{-1}]/(r'_1 r'_{-1} = w)
= \CC[r'_1,r'_{-1}]$. Thus $\cM_C$ is $\CC^2 = T^*\CC$.

(2)' If we replace $\bM$ by $\CC^{2\ell} = T^*(\CC^\ell)$, we get
\begin{equation*}
    \Spec H_*^{\bG(\cO)}(\cR_{\bG,\bN}) = \CC[w,r'_1,r'_{-1}]/ (r'_1 r'_{-1} = w^\ell).
\end{equation*}
Thus $\cM_C$ is type $A_{\ell-1}$ simple singularity.

(2)'' If we replace the $\bG$ action on $\bN$ by the weight $\ell$ representation, we get
$\Spec H_*^{\bG(\cO)}(\cR_{\bG,\bN}) = \CC[w,r'_1,r'_{-1}]/ (r'_1 r'_{-1} = w^\ell)$.
Thus $\cM_C$ is type $A_{\ell-1}$ simple singularity.

(3) Suppose $\bG = \CC^\times$, $\bM = T^*\CC^\times$ ($\bN=\CC^\times$). $\bG$ acts on $\bN$ by
the multiplication. Then 
\begin{equation*}
    \cR_{\bG,\bN} = \bigsqcup_{n\in\ZZ} \{ \left([z^n],s'(z)\right)
    \in\Gr_\bG\times \bN(\cO)\mid z^{-n}s'(z)\in\bN(\cO) \},
\end{equation*}
Since $\bN(\cO) = \{ s(z)\in\CC[[z]]\setminus \{0\} \mid 1/s(z)\in\CC[[z]]\} = \{ s(z) = s_0 + s_1 z + \dots\mid s_0\neq 0\}$, the component for $n\neq 0$ is $\emptyset$. The component for
$n=0$ is $\bN(\cO)$. Therefore
\begin{equation*}
    H_*^{\bG(\cO)}(\cR_{\bG,\bN})\cong
    H_*^{\bG(\cO)}(\bN(\cO)) = H_*(\mathrm{pt}) = \CC.
\end{equation*}
Thus $\cM_C = \mathrm{pt}$.
\end{Example}

As we see above, $\cR_{\bG,\bN}$ could have many connected components. In fact, we
have $\pi_0(\cR_{\bG,\bN}) = \pi_0(\Gr_\bG) \cong \pi_1(\bG)$. In above examples,
$\pi_1(\CC^\times) = \ZZ$. Since $\cR_{\bG,\bN}$ decomposes according to $\pi_1(\bG)$,
$\cM_C = \Spec H^{\bG(\cO)}_*(\cR_{\bG,\bN})$ is 
equipped with the action of $\pi_1(\bG)^\wedge$, the Pontryagin dual of $\pi_1(\bG)$.

\subsection{Ring objects}\label{subsec:ring}

Let $D_{\bG(\cO)}(\Gr_\bG)$ be the derived category of $\bG(\cO)$-equivariant constructible
sheaves on $\Gr_\bG$, more precisely its ind-completion.
It is equipped with a natural monoidal convolution structure $\ast$.
If $\omega_{\cR_{\bG,\bN}}$ is the dualizing sheaf of $\cR_{\bG,\bN}$,
the pushforward $\cA \defeq \pi_*(\omega_{\cR_{\bG,\bN}})$ is an object in $D_{\bG(\cO)}(\Gr_\bG)$,
and $H^{\bG(\cO)}_*(\cR_{\bG,\bN})$ is naturally isomorphic to
$H^*_{\bG(\cO)}(\Gr_\bG, \cA)$. Moreover, $\cA$ is a commutative ring object in $D_{\bG(\cO)}(\Gr_\bG)$,
i.e., it is equipped with a homomorphism $\cA\ast\cA\to\cA$ in $D_{\bG(\cO)}(\Gr_\bG)$.
It is `commutative' in the sense that it is compatible with the factorization given by Beilinson-Drinfeld Grassmanian.

\begin{Lemma}\label{lem:ring_objects}
    Suppose $\varphi\colon \bG_1\to\bG_2$ is a group homomorphism and
    $\Gr\varphi\colon\Gr_{\bG_1}\to\Gr_{\bG_2}$ be the corresponding morphism
    between affine Grassmannians.

    \textup{(1)} If $\cA_2$ is a commutative ring object in $D_{\bG_2(\cO)}(\Gr_{\bG_2})$,
    $\Gr\varphi^!\cA_2$ is a commutative ring object in $D_{\bG_1(\cO)}(\Gr_{\bG_1})$.

    \textup{(2)} If $\cA_1$ is a commutative ring object in $D_{\bG_1(\cO)}(\Gr_{\bG_1})$,
    $\Gr\varphi_*\cA_1$ is a commutative ring object in $D_{\bG_2(\cO)}(\Gr_{\bG_2})$.
\end{Lemma}

In particular, if $\Delta\bG\to\bG\times\bG$ is the diagonal embedding, commutative ring objects
$\cA_1$, $\cA_2$ give rise another commutative ring object $\cA_1\otimes^!\cA_2 = \Gr\Delta^!(\cA_1\boxtimes\cA_2)$.

If $\varphi\colon \bG\to\{e\}$ is the homomorphism to the trivial group, $\Gr\varphi_*(\cA)$ is a commutative
ring object in $D(\Gr_{\{e\}})$. But this derived category is the category of graded vector spaces, hence
$\Gr\varphi_*(\cA)=H^*_{\bG(\cO)}(\Gr_\bG,\cA)$ is a commutative ring over $\CC$ equipped with a grading compatible with multiplication.
Hence we can define an affine scheme by $\Spec H^*_{\bG(\cO)}(\Gr_\bG,\cA)$. This construction generalizes the definition of Coulomb branches.

\subsection{Regular sheaf}\label{subsec:reg}

Let $\mathrm{Perv}_{\bG(\cO)}(\Gr_\bG)$ be the
category of $\bG(\cO)$-equivariant perverse sheaves
on $\Gr_\bG$. It is an abelian category equipped with
a tensor structure given by the convolution $\ast$.
It is the heart of a $t$-structure on $D_{\bG(\cO)}(\Gr_\bG)$.
Geometric Satake equivalence says that there is an
equivalence of tensor categories
\begin{equation*}
    (\mathrm{Perv}_{\bG(\cO)}(\Gr_\bG),\ast)
    \cong
    (\operatorname{\mathscr Rep} \bG^\vee, \otimes),
\end{equation*}
where the right hand side is the tensor category
of finite dimensional representations of the
Langlands dual group $\bG^\vee$. We take ind-completion
in both sides, and consider $\CC[\bG^\vee]$ the
regular representations of $\bG^\vee$ in the right
hand side. Let $\cA_{\CC[\bG^\vee]}$ be the
corresponding object in the left hand side.
We call it the \emph{regular sheaf}.
It is a ring object in the left hand side, and 
is equipped with a $\bG^\vee$ action, induced
from the \emph{right} multiplication action on $\CC[\bG^\vee]$.
Note that the left $\bG^\vee$ action is used
in the geometric Satake equivalence, but
the right action remains.
More concretely
\begin{equation*}
    \cA_{\CC[\bG^\vee]}\cong\bigoplus
    V_{G^\vee}(\lambda)^*\otimes\operatorname{IC}^\lambda,
\end{equation*}
where $V_{G^\vee}(\lambda)$ is the irreducible
represenation of $G^\vee$ with highest weight $\lambda$,
and
$\operatorname{IC}^\lambda$ is the intersection
cohomology complex for the orbit
$\bG(\cO)\cdot z^\lambda$ with $\lambda$ viewed
as a coweight of $\bG$. The summation runs over
the set of dominant integral coweights of $\bG$.

For example, $\cA_{\CC[\bG^\vee]}$ is the
constant sheaf if $\bG=\CC^\times$, and
is the direct sum of constant sheaves on $\{[z^n]\}$
($n\in\ZZ$).

\begin{Theorem}[\cite{MR2053952}]\label{thm:3}
    \textup{(1)} Let $\iota\colon\{1\}=\Gr_{\{e\}}\to
    \Gr_{\bG}$ be the inclusion induced from $\{1\}\to\bG$. Then $\Spec\iota^!\cA_{\CC[\bG^\vee]}$
    is isomorphic to the nilpotent
    cone $\CN_{\bG^\vee}$ of $\bG^\vee$.

    \textup{(2)} $\Spec H^*_{\bG(\cO)}(\Gr_{\bG},
    \cA_{\CC[\bG^\vee]})$ is isomorphic to
    $\bG^\vee\times\mathcal S^\vee$, where $\mathcal S^\vee$
    is the Kostant-Slodowy slice to a principal
    nilpotent orbit of $\bG^\vee$.
\end{Theorem}

Recall that we take an $\algsl_2$-triple
$e$, $f$, $h$ for
a nilpotent element
$e\in\CN_{\bG^\vee}$ so that
$[h,e] = 2e$, $[h,f] = - 2f$, $[e,f] = h$, and define a \emph{slice}
by the affine space
\begin{equation*}
    e + \Ker \operatorname{ad}(f) \subset
    \Lie\bG^\vee.
\end{equation*}
The above $\mathcal S^\vee$ is defined for a principal
nilpotent element $e$. 
(In fact, it is more natural
to define the slice in the \emph{dual} of $\Lie\bG^\vee$,
as $\bG^\vee\times\mathcal S^\vee$ is a symplectic reduction
of $T^*\bG^\vee$.)
The inclusion $\mathcal S^\vee\hookrightarrow\Lie\bG^\vee$
induces an isomorphism $\mathcal S^\vee\cong \Lie\bG^\vee/
\operatorname{Ad}\bG^\vee = \Lie\mathbf T^\vee/\mathbb W$,
where $\mathbf T^\vee$ is a maximal torus of $\bG^\vee$
and $\mathbb W$ is the Weyl group of $\bG^\vee$.
As we remarked above, $\cA_{\CC[\bG^\vee]}$ has a $\bG^\vee$ action which
is identified with natural $\bG^\vee$ actions on $\CN_{\bG^\vee}$
and $\bG^\vee\times\mathcal S^\vee$ above respectively.

We have two ways to construct ring objects. It is natural to ask whether there is
a way to realize $\cA_{\CC[\bG^\vee]}$ from $\bG\curvearrowright\bM$. Here is an example. Let
\begin{gather*}
    \bN \defeq \Hom(\CC,\CC^2)\oplus \Hom(\CC^2,\CC^3)\oplus\cdots\oplus
    \Hom(\CC^{n-1},\CC^n),\\
    \bG = \operatorname{PGL}_n, \qquad \tilde{\bG} = 
    \left(\GL_1\times\GL_2\times\cdots\times
    \GL_n\right)/\CC^\times,
\end{gather*}
where $\CC^\times$ is the diagonal central subgroup. We have a ring homomorphism
$\varphi\colon\tilde{\bG}\to\bG$. Let $\cA_{\tilde{\bG},\bM}$ be the pushforward
of the dualizing sheaf on $\cR_{\tilde{\bG},\bN}$ associated with $(\tilde{\bG},\bN)$.
It is a ring object in $D_{\tilde{\bG}(\cO)}(\Gr_{\tilde\bG})$. We then pushforward it
to $\Gr_\bG$, namely $\Gr\varphi_*\cA_{\tilde{\bG},\bM}$. Note that $\bG^\vee=\SL_n$.

\begin{Theorem}[\protect{\cite[Th.~2.11]{2017arXiv170602112B}}]
\label{thm:ringobj}
    $\Gr\varphi_*\cA_{\tilde{\bG},\bM}$ and $\cA_{\CC[\SL_n]}$ are isomorphic
    as ring objects in $D_{\operatorname{PGL}_n(\cO)}(\Gr_{\PGL_n})$.
\end{Theorem}

It is expected that we have similar results for $\bG=\SO_{2n}$, $\grpSp_{2n}$
in view of \cite[\S5.2]{MR2610576}. However, there is no analogous result for
other groups, including $\bG=\SO_{2n+1}$.

\subsection{S-dual}

Let us now give the (provisual) definition of the S-sual of $\bG\curvearrowright\bM$ under
the assumption $\bM = T^*\bN$. We first form the ring object given by the pushforward of the dualizing sheaf of $\cR_{\bG,\bN}$ as in \S\ref{subsec:ring}. Let us denote it by
$\cA_{\bG,\bM}$.
We define a ring object $\cA_{\bG,\bM}\otimes^!\cA_{\CC[\bG^\vee]}$. Then
we consider
\begin{equation}\label{eq:1}
    \bM^\vee \defeq \Spec H^*_{\bG(\cO)}(\Gr_\bG, \cA_{\bG,\bM}\otimes^!\cA_{\CC[\bG^\vee]}).
\end{equation}
It is equipped with $\bG^\vee$ action from $\cA_{\CC[\bG^\vee]}$. Since the S-dual
is defined for the pair $(\bG,\bM)$, we say $(\bG^\vee\curvearrowright\bM^\vee)$ is
the S-dual of $(\bG\curvearrowright\bM)$.
The definition \eqref{eq:1} is available only when
$\bM=T^*\bN$, but we expect that it can be weakened
to more general $\bM$ as we mentioned in \S\ref{sec:intro}.

\begin{Remarks}\label{rem:S-dual_phys}
    (1) The definition \eqref{eq:1} is derived from 
    $\cT_{\bG,\bM}^\vee = \left(\left(\cT_{\bG,\bM}\times\cT[\bG]\right)
    \tslabar\bG\right)^*$
    (\cite[(3.13)]{MR2610576}) in a physics level of rigor as
    follows. If this would be a theory defined by
    $\bG^\vee\curvearrowright\bM^\vee$, its Higgs branch
    is $\bM^\vee$. Since $*$ is the 3d mirror symmetry, exchanging
    Higgs and Coulomb branches, the Coulomb branch of
    $\left(\cT_{\bG,\bM}\times\cT[\bG]\right)
    \tslabar\bG$ should be $\bM^\vee$. It is $H^*_{\bG(\cO)}(\Gr_\bG,\cA_{\cT_{\bG,\bM}\times\cT[\bG]})$ of the ring object associated with $\cT_{\bG,\bM}\times\cT[\bG]$, which is the $!$ tensor product
    of $\cA_{\bG,\bM}$ and $\cA_{\cT[\bG]}$. Then
    $\cA_{\cT[\bG]}\cong \cA_{\CC[\bG^\vee]}$ is a proposal in \cite{2017arXiv170602112B}.

    (2) As in Remark~\ref{rem:just}, we need a justification of
    the notation $\cA_{\bG,\bM}$, as we use the description
    $\bM = T^*\bN$. If $\bN$ is a representation of $\bG$ as in
    \cite{bdfrt}, $\cA_{\bG,\bM}$ is defined without using
    the description $\bM = \bN\oplus\bN^*$, but it yet depends on the choice of a square root $\sqrt{\mathscr{L}}$
    of the determinant line bundle $\mathscr{L}$. (See
    \cite[\S4.2]{bdfrt} for an identification with the construction
    from $\cR_{\bG,\bN}$.)
    Therefore it is denoted by $\cA_{\bG,\bM,\sqrt{\mathscr{L}}}$.
\end{Remarks}

\begin{Example}\label{ex:2}
\textup{(1)} Suppose $\bG$ is a torus. Then $\cA_{\CC[\bG^\vee]}$ is the constant sheaf
over $\Gr_{\bG} = \Hom(\CC^\times,\bG)$, which is also the dualizing sheaf. Therefore
$\cA_{\bG,\bM}\otimes^!\cA_{\CC[\bG^\vee]} = \cA_{\bG,\bM}$. Therefore $\bM^\vee$ is
nothing but the Coulomb branch $\cM_C(\bG,\bM)$. The action of $\bG^\vee$ is identified
with the action of $\pi_1(\bG)^\wedge$.

\textup{(2)} Theorem~\ref{thm:3}(2) can be understood as the S-dual of 
$(\bG\curvearrowright\mathrm{pt})$ is $(\bG^\vee\curvearrowright\bG^\vee\times\mathcal S^\vee)$.
Since the S-dual of the S-dual is expected to be the original space, we conjecture that
the S-dual of $(\bG^\vee\curvearrowright\bG^\vee\times\mathcal S^\vee)$ is $(\bG\curvearrowright\mathrm{pt})$. If $\bG=\CC^\times$, Example~\ref{ex:1}(3) confirms
that it is true.

\textup{(3)} Consider $\bG\curvearrowright \bM = T^*\bG$, where $\bG$ acts on $\bN=\bG$ by the left multiplication. Then
$\cR_{\bG,\bN} = \{ [g(z),s(z)]\in \bG(F)\times^{\bG(\cO)}\bG(\cO)
\mid g(z)s(z)\in\bG(\cO)\}$ is nothing but $\bG(\cO)$, as
we have $g(z)\in\bG(\cO)$. Furthermore, $\cR_{\bG,\bN}$ is
mapped to $1\in\Gr_{\bG}$. Therefore $\cA_{\bG,\bM}\otimes^!\cA_{\CC[\bG^\vee]}$ is 
$H^*_{\bG(\cO)}(\bG(\cO),\omega_{\bG(\cO)})\otimes
\iota^!\cA_{\CC[\bG^\vee]}$, where $\iota\colon\{1\}\to\Gr_{\bG}$
is the inclusion. By Theorem~\ref{thm:3}, this is isomorphic to
$\CC[\CN_{\bG^\vee}]$, the coordinate ring of the nilpotent cone
of $\bG^\vee$. Thus the S-dual is $(\bG^\vee\curvearrowright\CN_{\bG^\vee})$. Note that
$\CN_{\bG^\vee}$ is singular unless $\bG^\vee$ is torus.

\textup{(4)} Since we expect the S-dual of the S-dual is the original space, we guess that the S-dual of $(\bG^\vee\curvearrowright\CN_{\bG^\vee})$ is $(\bG\curvearrowright T^*\bG)$. Since $\CN_{\bG^\vee}$ is singular (unless $\bG^\vee$
is torus), we cannot check this statement. Let us speculate
why this could be true.

By \cite[5(x)]{2017arXiv170602112B} we know that
$H^*_{\bG^\vee(\cO)}(\Gr_{\bG^\vee}, \cA_{\CC[\bG]}\otimes^!
\cA_{\CC[\bG]})\cong \CC[T^*\bG]$. Therefore the statement
can be deduced \emph{formally}, if we `show' that
the ring object
$\cA$ attached to $\bG^\vee\curvearrowright\CN_{\bG^\vee}$ is $\cA_{\CC[\bG]}$.
Since $\CN_{\bG^\vee}$ is not smooth, not the cotangent bundle,
the current technology does not give a definition of $\cA$
for  $\bG^\vee\curvearrowright\CN_{\bG^\vee}$.
We will return back this claim in Example~\ref{ex:3}.
\end{Example}

Let us explain the relation between \eqref{eq:1} and the derived 
Satake equivalence. The formula \eqref{eq:1}, without
taking $\Spec$, i.e.,
$H^*_{\bG(\cO)}(\Gr_\bG,\mathcal F\otimes^!\cA_{\CC[\bG^\vee]})$
makes sense for an object $\mathcal F$ of
$D_{\bG(\cO)}(\Gr_\bG)$. (The result is not necessarily a ring,
hence we cannot take $\Spec$.)
In \cite[Lemma~5.13]{2017arXiv170602112B} we showed
that the functor is the inverse of the derived Satake equivalence in \cite{MR2422266} up to the twist
by the Chevalley involution. See 
\cite[\S5(vi),\S5(vii)]{2017arXiv170602112B} for 
the framework and the precise statement. 

\section{Symplectic reduction}

\subsection{Definition}

Suppose that a product $\bG_1\times\bG_2$ (resp.\ $\bG_2\times\bG_3$) acts on a smooth affine symplectic variety
$\bM_{12}$ (resp.\ $\bM_{23}$). 
To simplify the statement below, we twist the $\bG_2$
(resp.\ $\bG_3$)
action on $\bM_{12}$ (resp.\ $\bM_{23}$) through an involution
on $\bG_2$ (resp.\ $\bG_3$) interchanging conjugacy clases of $g$ and
$g^{-1}$ (the Chevalley involution).
We symbolically denote the space with action as
$\bG_1\curvearrowright\bM_{12}\curvearrowleft\bG_2$
(resp.\ $\bG_2\curvearrowright\bM_{23}\curvearrowleft\bG_3$).
We consider the symplectic reduction
\begin{equation*}
    \bM_{12}\circ\bM_{23} \defeq
    \mu^{-1}_2(0)\dslash \bG_2,
\end{equation*}
where $\mu_2$ is the moment map for the $\bG_2$-action
on the product $\bM_{12}\times\bM_{23}$,
and $\dslash\bG_2$ is the categorical quotient. This is
an affine algebraic variety, possibly with singularities.
If $\bG_2$ action is free, the reduction is smooth, and is equipped with a symplectic form induced from those of $\bM_{12}$
and $\bM_{23}$.
Therefore we regard $\bM_{12}\circ\bM_{23}$ as a symplectic manifold
in a generalized sense without any justification.
Moreover, since we took the quotient by $\bG_2$, we have
the induced action $\bG_1\curvearrowright\bM_{12}\circ\bM_{23}
\curvearrowleft\bG_3$.
Thus we can define the `composition' $(\bM_{12}\circ\bM_{23})\circ\bM_{34}$ 
of the reduction as $\bM_{12}\circ\bM_{23}\circ\bM_{34}$, 
even if $\bM_{12}\circ\bM_{23}$ may \emph{not} be smooth.

%
%
%
If $\bM_{12} = T^*\bN_{12}$, $\bM_{23} = T^*\bN_{23}$,
and the $\bG_2$-action on $\bN_{12}\times\bN_{23}$ is
free, we have $\bM_{12}\circ \bM_{23} = T^*\left(\frac{\bN_{12}\times\bN_{23}}{\bG_2}\right)$.

\subsection{Ring object associated with the reduction}\label{subsec:ring_red}

Suppose that a product $\bG_1\times\bG_2$ (resp.\ $\bG_2\times\bG_3$) acts on a smooth affine variety
$\bN_{12}$ (resp.\ $\bN_{23}$). 
%
%
We 
denote the space with action by 
$\bG_1\curvearrowright\bN_{12}\curvearrowleft\bG_2$
(resp.\ $\bG_2\curvearrowright\bN_{23}\curvearrowleft\bG_3$)
as before.
We consider the ring object $\cA_{12} = \cA_{\bG_1\times\bG_2,\bM_{12}}$
(resp.\ $\cA_{23} = \cA_{\bG_2\times\bG_3,\bM_{23}}$) on
$\Gr_{\bG_1\times\bG_2}$ (resp.\ $\Gr_{\bG_2\times\bG_3}$).
We define a ring object on $\Gr_{\bG_1\times\bG_3}$ by
\begin{equation*}
    \cA_{12}\circ\cA_{23}
    \defeq p_{13*}\left(p_{12}^!\cA_{12}\otimes^!
    p_{23}^!\cA_{23}\right),
\end{equation*}
where $p_{ab}\colon \Gr_{\bG_1\times\bG_2\times\bG_3}
\cong\Gr_{\bG_1}\times\Gr_{\bG_2}\times\Gr_{\bG_3}
\to \Gr_{\bG_a}\times\Gr_{\bG_b}$ is the projection
to the product of the $a$th and $b$th factors.
If the action of $\bG_2$ on $\bN_{12}\times\bN_{23}$
is free, $\cA_{12}\circ\cA_{23}$ is the ring object
associated with the quotient space
$
\frac{\bN_{12}\times\bN_{23}}{\bG_2}$.
Recall that $T^*\left(\frac{\bN_{12}\times\bN_{23}}{\bG_2}\right)$ is the symplectic reduction
$\bM_{12}\circ\bM_{23}$.
%
Thus we regard $\cA_{12}\circ\cA_{23}$ as the ring object
associated with $\bG_1\curvearrowright(\bM_{12}\circ\bM_{23})
\curvearrowleft\bG_3$, even when
$\bM_{12}\circ\bM_{23}$ is singular.
In the terminology of the gauge theory, $\cA_{12}\circ\cA_{23}$
is the (Coulomb) branch ring object of the 3d $\mathcal N=4$ 
theory $\left(\mathcal T_{\bG_1\curvearrowright\bM_{12}\curvearrowleft\bG_2}
\times
\mathcal T_{\bG_2\curvearrowright\bM_{23}\curvearrowleft\bG_3}
\right)
\tslabar\bG_2$, obtained from
the product of theories associated with $\bM_{12}$ and 
$\bM_{23}$ by the gauging by $\bG_2$.

\begin{Example}\label{ex:3}
    Consider 
    \begin{gather*}
     \bM_0\circ\bM_1\circ\dots\circ\bM_{n-1}\\
     \text{ where } \GL_i\curvearrowright\bM_i = T^*\Hom(\CC^i,\CC^{i+1})
    \curvearrowleft\GL_{i+1}
    \end{gather*}
    By \cite{MR549399},
    $\bM_0\circ\bM_1\circ\dots\circ\bM_{n-1}\cong\CN_{\GL_n}$
    as $\GL_0\times\GL_n = \GL_n$-varieties.
    By $\GL$-version of Theorem~\ref{thm:ringobj}, the ring
    object corresponding to $\bM_0\circ\bM_1\circ\dots\circ\bM_{n-1}$
    in the above sense is $\cA_{\CC[\GL_n]}$. This confirms
    the prediction in Example~\ref{ex:2}(4) for $\bG = \GL_n$.
\end{Example}

\subsection{S-duality commutes with reduction}

Suppose we have
$\bG_1\curvearrowright\bM_{12} = T^*\bN_{12}\curvearrowleft\bG_2$
and $\bG_2\curvearrowright\bM_{23} = T^*\bN_{23}\curvearrowleft\bG_3$
as before.
We have ring objects $\cA_{12}$, $\cA_{23}$, and we can define
S-dual by \eqref{eq:1}:
$\bG_1^\vee\curvearrowright\bM_{12}^\vee = \Spec H^*_{(\bG_1\times\bG_2)(\cO)}(
\Gr_{\bG_1\times\bG_2},\cA_{\CC[\bG_1^\vee]}\otimes^!
\cA_{12}\otimes^!\cA_{\CC[\bG_2^\vee]})\curvearrowleft\bG_2^\vee
$
and similarly $\bG_2^\vee\curvearrowright\bM_{23}^\vee
\curvearrowleft\bG_3^\vee$.

The following result says that the S-duality commutes with
the reduction.

\begin{Theorem}[\protect{\cite[Remark~5.21]{2017arXiv170602112B}}]
\label{thm:B}
We have a $\bG_1^\vee\times\bG_3^\vee$-equivariant isomorphism
\begin{equation*}
    (\bM_{12}\circ\bM_{23})^\vee \cong \bM_{12}^\vee\circ\bM_{23}^\vee,
\end{equation*}
where the left hand side is understood as
\begin{equation*}
    \Spec H^*_{(\bG_1\times\bG_3)(\cO)}(\Gr_{\bG_1\times\bG_3}, \cA_{\CC[\bG_1^\vee]}\otimes^!
    (\cA_{12}\circ\cA_{23})\otimes^!\cA_{\CC[\bG_3^\vee]})
\end{equation*}
by \S\ref{subsec:ring_red}.
\end{Theorem}

Recall that 
we regard $\cA_{12}\circ\cA_{23}$ as the ring object
associated with $\bM_{12}\circ\bM_{23}$, 
even though $\bM_{12}\circ\bM_{23}$ is singular in general.

\begin{Remark}
    This result is manifest for quantum field theories obtained
    by configurations of D3,D5/NS5 branes \cite{MR1451054}, which underlie results in the next section. 
    (D5 is $\bM_\times$, NS5 is $\bM_\circ$.)
    The operation $\circ$ corresponds
    to the concatenation of two configurations. The S-duality
    exchanges D5 and NS5 branes. The result simply means that
    concatenation and the exchange D5 $\leftrightarrow$ NS5 commute. 
\end{Remark}

\section{Bow varieties and Coulomb branches}

In \cite{2016arXiv160602002N}, we proved that
Coulomb branches of quiver gauge theories
of affine type $A$ are Cherkis bow varieties. Cherkis bow varieties are products of two types of
algebraic $\GL(V_i)\times\GL(V_j)$-Hamiltonian manifolds:
\begin{equation*}
    \begin{split}
        &\GL(V_i)\curvearrowright
        \bM_{\medcirc}(V_i,V_j) =
        T^*\Hom(V_i,V_j) \curvearrowleft\GL(V_j)\\
        &\GL(V_i)\curvearrowright
        \bM_{\times}(V_i,V_j) =
        \begin{cases}
        \GL(V_i)\times\mathcal S(\bv_i-\bv_j,1^{\bv_j})
            & \text{if $\bv_i > \bv_j$}\\
        (i\longleftrightarrow j) & \text{if $\bv_i < \bv_j$} \\
        T^*(\GL(V_i)\times V_i) & \text{if $\bv_i=\bv_j$}
        \end{cases} 
        \curvearrowleft\GL(V_j).
    \end{split}
\end{equation*}
Let us explain the notation appearing in the
definition of $\bM_\times(V_i,V_j)$. We
set $\bv_i = \dim V_i$, $\bv_j = \dim V_j$.
If $\bv_i > \bv_j$, 
$\mathcal S(\bv_i-\bv_j,1^{\bv_j})$ is the 
slice to the nilpotent orbit corresponding
to the hook partition $(\bv_i-\bv_j,1^{\bv_j})$
of $\bv_i$, as in \S\ref{subsec:reg}.
The action of $\GL(V_i)$ is the right multiplication to
the factor $\GL(V_i)$ in $\bM_\times(V_i,V_j)$.
The centralizer of the corresponding $\algsl_2$ triple
is $\GL(V_j)$, acting on $\mathcal S(\bv_i-\bv_j,1^{\bv_j})$ naturally.
We let it act on the factor $\GL(V_i)$ by the left multiplication.
If $\bv_i < \bv_j$, we exchange the roles of $i$, $i$.
(Strictly speaking, we used a variant of $\mathcal S(\bv_i-\bv_j,1^{\bv_j})$ was used in \cite{2016arXiv160602002N}.
But they give an isomorphic hamiltonian space.)

The identification of Coulomb branches with bow varieties were proved without reference to the S-dual, but it follows from
the following:
\begin{Conjecture}\label{conj:S-dual_bow}
    Two hamiltonian spaces $\GL(V_i)\curvearrowright
        \bM_{\medcirc}(V_i,V_j)
        \curvearrowleft\GL(V_j)$ and
        $\GL(V_i)\curvearrowright
        \bM_{\times}(V_i,V_j)
        \curvearrowleft\GL(V_j)$
        are S-dual to each other.
\end{Conjecture}

A quiver gauge theory of finite type $A_\ell$ is given by the
product
$\bM_{\medcirc}(V_i,V_{i+1})$ ($i=1,2,\dots,\ell-1$),
together with $\bM_{\medcirc}(V_i,W_i)$ ($i=1,\dots,\ell$) 
for a `framing vector space' $W_i$,
gauged by $\prod_i \GL(V_i)$:
\begin{equation*}
     \left(\prod{\bM_{\medcirc}(V_i,V_{i+1})}
     \oplus{\bM_{\medcirc}(V_i,W_i)}\right)
     \tslabar \prod_i \GL(V_i).
\end{equation*}
Here we confuse hamiltonian spaces and the corresponding
quantum field theories bu abuse of notation.
The factor $\bM_{\medcirc}(V_i,W_i)$ can be replaced by
$\bM_{\times}(V_i,V_i)\times\cdots\times \bM_{\times}(V_i,V_i)$
($\bw_i = \dim W_i$ factors). 
Therefore it is of the form
\begin{equation}\label{eq:6}
    \bM_\medcirc(0,V_1)\circ
    \underbrace{\bM_\times(V_1,V_1)\circ\dots
    \circ\bM_\times(V_1,V_1)}_{\text{$\bw_1$ times}}\circ
    \bM_\medcirc(V_1,V_2)\circ\cdots.
\end{equation}

Therefore the corresponding Coulomb branch
is given by the successive application of Theorem~\ref{thm:B}.
By the conjecture above, we get
\begin{equation}\label{eq:5}
    \bM_\times(0,V_1)\circ
    \underbrace{\bM_\medcirc(V_1,V_1)\circ\dots
    \circ\bM_\times(V_1,V_1)}_{\text{$\bw_1$ times}}\circ
    \bM_\times(V_1,V_2)\circ\cdots.
\end{equation}
This is, indeed, the main result of \cite{2016arXiv160602002N}.

The proof in \cite{2016arXiv160602002N} does not use
Conjecture~\ref{conj:S-dual_bow}: we directly show that
\eqref{eq:5} is the Coulomb branch of the
gauge theory \eqref{eq:6}. In turn, we can show
Conjecture~\ref{conj:S-dual_bow} in
the cases $\bM_\times(V_i,V_j)$ with $\bv_i = \bv_j$ and $\bM_\circ(V_i,V_j)$ with arbitrary $\bv_i$, $\bv_j$
from this statement together
with the realization of $\cA_{\CC[\GL_n]}$ by Theorem~\ref{thm:ringobj}.
Indeed, we add $\cdots\circ
\bM_\medcirc(\CC^{\bv_i-2},\CC^{\bv_i-1})
\circ\bM_\medcirc(\CC^{\bv_i-1},V_i)\circ$ and
$\circ\bM_\medcirc(V_j,\CC^{\bv_j-1})\circ
\bM_\medcirc(\CC^{\bv_j-1},\CC^{\bv_j-2})\circ\cdots$
to the left and right of $\bM_\medcirc(V_i,V_j)$. 
The Coulomb branch of this theory is the
S-dual of $\GL(V_i)\curvearrowright\bM_\medcirc(V_i,V_j)
\curvearrowleft\GL(V_j)$. 
The theory is of the form of \eqref{eq:6}, hence the Coulomb branch
is given by \eqref{eq:5}. Since the added parts give
$T^*\GL(V_i)$ and $T^*\GL(V_j)$ respectively, the answer 
\eqref{eq:5} is nothing but $\bM_\times(V_i,V_j)$ or
$\bM_\medcirc(V_i,V_j)$ with $\bv_i=\bv_j$.

In the remaining case $\bM_\times(V_i,V_j)$ with $\bv_i\neq\bv_j$,
we add $\cdots\circ
\bM_\medcirc(\CC^{\bv_i-2},\CC^{\bv_i-1})
\circ\bM_\medcirc(\CC^{\bv_i-1},V_i)\circ$ and
$\circ\bM_\medcirc(V_j,\CC^{\bv_j-1})\circ
\bM_\medcirc(\CC^{\bv_j-1},\CC^{\bv_j-2})\circ\cdots$ as above,
and use the Hanany-Witten transition (see \cite[\S7.2]{2016arXiv160602002N}) to a gauge theory of the form \eqref{eq:6}.
Then we apply the Hanany-Witten transition to \eqref{eq:5}.
Since Hanany-Witten transition
is compatible with the exchange of $\times$ and $\medcirc$,
we get $\bM_\circ(V_i,V_j)$.
Thus we only need to check that the ring object $\cA$
is independent under Hanany-Witten transition.

\section{Hyperspherical conditions}\label{sec:hyper}

In Remark~\ref{rem:S-dual_phys}(1), we said that
$\bG^\vee\curvearrowright\bM^\vee$ is the Higgs
branch of the theory
$\cT_{\bG,\bM}^\vee = \left(\left(\cT_{\bG,\bM}\times\cT[\bG]\right)
    \tslabar\bG\right)^*$. 
If this theory comes from a \emph{smooth} $\bM^\vee$,
the corresponding Coulomb branch must be a single point. Thus
the Higgs branch of $\left(\cT_{\bG,\bM}\times\cT[\bG]\right)
    \tslabar\bG$ must be a point, i.e.,
\begin{equation*}
    \left(\bM\times\mathcal N_\bG\right)\tslash\bG
\end{equation*}
is a point.
 
One of hyperspherical conditions says
\begin{itemize}
    \item a generic $\bG$ orbit in $\bM$ is a coisotropic subvariety.
\end{itemize}
See \cite[\S3.5]{BZSV}. 
This is equivalent to the condition
\begin{itemize}
    \item All nonempty fibers of 
    $\mu\dslash\bG\colon \bM\dslash\bG\to
    \operatorname{Lie}\bG^*\dslash\bG$ are
    finite
\end{itemize}
by \cite{MR2545864}.\footnote{The author thanks
Ivan Losev for his explanation how to derive the
equivalence from \cite{MR2545864}.}

\bibliographystyle{amsalpha}
\bibliography{nakajima,mybib,orthsymp,kika}

\section*{Errata}

\begin{itemize}
    \item 1.~Introduction, the second paragraph, l.7 : replace 'pairs' by `the pair'.
    \item 1.~Introduction, the paragraph after \eqref{eq:7} : regarding \cite{bdfrt,2022arXiv220901088T}, we should impose the \emph{anomaly free} condition in order to have a well-defined Coulomb branch.
    \item Example~\ref{ex:Hitchin}, l.4 : Delete `is'.
    \item Example~\ref{ex:1} (2)', (2)'' : these two examples give
    $\bM^\vee\cong A_{\ell-1}\cong \bM^{\prime\vee}$ with $\bM = T^*\CC^{\ell}$, $\bM' = T^*\CC$ (weight $\ell$) in view of Example~\ref{ex:2}(1). They suggest that a finer structure of the S-dual $A_{\ell-1}$ is necessary in order to recover the original spaces $\bM$, $\bM^\vee$ from the S-dual.
    At this moment, it is not clear what is the finer structure. In these
    examples, note that (2)' realizes $A_{\ell-1}$ as the symplectic
    reduction of $\CC^{2\ell}$ by an $\ell-1$ dimensional torus, while
    (2)'' realizes it as the quotient $\CC^2/(\ZZ/\ell\ZZ)$ by
    \cite[\S3(vii)(d)]{2016arXiv160103586B} and \S3(vii)(c) 
    respectively.
    
    \item Lemma~\ref{lem:ring_objects}(2), we assume that $\varphi$ is surjective to make $\Gr\varphi_*\cA_1$ is
    $\bG_2(\cO)$-equivariant.
    \item 3.4~S-dual, the first line : replace 'provisual' by `provisional'.
    \item Remarks 3.7 : Add the following: (3) The Coulomb branch $\cM_C(\bG,\bM)$ is recovered from the S-dual $\bM^\vee$ by the Kostant reduction, i.e., $\cM_C(\bG,\bM)\cong\left(\bM^\vee\times(\bG^\vee\times\cS^\vee)\right)\tslash\bG^\vee$. This follows from the proof of \cite[Lemma~5.18]{2017arXiv170602112B}.

    \item line -5 above Conjecture~\ref{conj:S-dual_bow} : replace one of $i$ by $j$.

    \item As a consequence of Conjecture~\ref{conj:S-dual_bow} and Theorem~\ref{thm:B}, we see
    that S-dual of $\GL_n\times\mathcal S_\lambda$ is $\overline{\mathbb O(\lambda^t)}$, where $\mathcal S_\lambda$ is the slice for a nilpotent element, whose Jordan form is given by the partition $\lambda$, and
    $\mathbb O(\lambda^t)$ is the nilpotent orbit, whose Jordan form
    is the transpose partition $\lambda^t$. Finally $\overline{\mathbb O(\lambda^t)}$ is the closure of $\mathbb O(\lambda^t)$. Indeed,
    Conjecture~\ref{conj:S-dual_bow} and Theorem~\ref{thm:B} give
    the S-dual as a reduction of a product of $\bM_\circ(V_i,V_{i+1})$ for
    certain $V_i$. It is exactly the realization of 
    $\overline{\mathbb O(\lambda^t)}$ in \cite{MR549399}.
    As a generalization of Example~\ref{ex:3}, we also see
    that the dual of $\overline{\mathbb O(\lambda^t)}$
    is $\GL_n\times\mathcal S_\lambda$.

    We conjecture that the same result also holds for more general
    complex semisimple groups $G$ and a slice $S_e$ for a nilpotent element $e$. Here $\overline{\mathbb O(\lambda^t)}$ is replaced by the affinization of a cover of the nilpotent orbit $\mathbb O_{d_{\mathrm{BV}}(e)}$, Barbasch-Vogan dual of the orbit of $e$.
    The cover is described in \cite{Chacaltana:2012zy,2021arXiv210803453L}.
    
    \item line 6 below Conjecture~\ref{conj:S-dual_bow} : replace `bu' by `by'.
    \item \eqref{eq:5} : Replace $\bM_\times(V_1,V_1)$ by $\bM_\circ(V_1,V_1)$.

    \item Section~\ref{sec:hyper} is not completed. Considering
    the fiber over $0\dslash\bG$ of $\mu\dslash\bG$, we see that
    one of hyperspherical conditions says that 
    $\left(\bM\times\mathcal N_\bG\right)\tslash\bG$ consists of
    finite points, but not necessarily a single point. We should explore hyperspherical conditions more in view of this article.
\end{itemize}

\end{document}